\documentclass[11pt]{article}
\usepackage{amsfonts}
\usepackage{amsmath}
 \title{\textbf{$3$-quasi-Sasakian manifolds}}
 \author{Beniamino Cappelletti Montano, Antonio De Nicola, Giulia Dileo}
   \date{}
\usepackage{graphicx}
\usepackage{amsmath}
\newtheorem{theorem}{Theorem}[section]

\newtheorem{corollary}[theorem]{Corollary}

\newtheorem{definition}[theorem]{Definition}
\newtheorem{example}[theorem]{Example}

\newtheorem{lemma}[theorem]{Lemma}

\newtheorem{remark}[theorem]{Remark}

\newenvironment{proof}[1][Proof]{\textbf{#1.} }{\ \rule{0.5em}{0.5em}}
\begin{document}

\maketitle

\begin{abstract}
In the present paper we carry on a systematic study of
$3$-quasi-Sasakian manifolds. In particular we prove that the three
Reeb vector fields generate an involutive distribution determining a
canonical totally geodesic and Riemannian foliation. Locally, the
leaves of this foliation turn out to be Lie groups: either the
orthogonal group or an abelian one. We show that $3$-quasi-Sasakian
manifolds have a well-defined rank, obtaining a rank-based
classification. Furthermore, we prove a splitting theorem for
these manifolds assuming the integrability of one of
the almost product structures. Finally, we show that the vertical distribution
is a minimum of the corrected energy.
\end{abstract}

\section{Introduction}

Quasi-Sasakian manifolds were introduced and deeply studied by D.
E. Blair in \cite{blair0} in the attempt to unify Sasakian and
cosymplectic geometry. He defined a quasi-Sasakian structure on a
$(2n+1)$-dimensional manifold $M$ as a normal almost contact
metric structure $(\phi,\xi,\eta,g)$ whose fundamental $2$-form
$\Phi$ is closed. His work was aimed especially to explore the
geometric meaning of the rank of the $1$-form $\eta$. For
instance, in a cosymplectic manifold $\eta$ has rank $1$ and in a
Sasakian manifold it has maximal rank $2n+1$. Blair proved that
there are no quasi-Sasakian manifolds of even rank. So, let $2p+1$
be the rank of a quasi-Sasakian manifold; in this case, if the
determined almost product structure is integrable, the space is
locally the product of a Sasakian manifold where
$\eta\wedge(d\eta)^p\neq 0$ and a K\"{a}hler manifold whose
fundamental $2$-form is $\Phi-d\eta$ properly restricted.

Next, some other general results were found by S. Tanno
\cite{tanno}, S. Kanemaki \cite{kanemaki1}, Z. Olszak
\cite{olszak2} and other authors. More recently quasi-Sasakian
manifolds appeared also in other contexts such as CR-geometry and
mathematical physics (see e.g. \cite{calin, friedrich, gonzales, rustanov}).
A study of quasi-Sasakian structures
using the techniques of $G$-structures by V. F. Kirichenko and A.
R. Rustanov should also be mentioned \cite{kirichenko}.

In the same years when the theory of quasi-Sasakian structures was
developed, the Japanese school including geometers such as Y. Y. Kuo, K. Yano, M.
Konishi, S. Ishihara et al. introduced the notion of almost
$3$-contact manifold as a triple of almost contact structures on a
$(4n+3)$-dimensional manifold $M$ satisfying the proper
generalization of the quaternionic identities.
From then on, several authors in the last 30--40 years have dealt
with almost $3$-contact geometry, in particular in the setting of
$3$-Sasakian manifolds, due to the increasing awareness of their
importance in mathematics and physics, together with the closely
linked hyper-K\"{a}hlerian and quaternion-K\"{a}hlerian manifolds.
 We refer the reader to the remarkable survey \cite{galicki2} and
references therein.

When each of the three structures of an almost $3$-contact metric
manifold is quasi-Sasakian, we say that the manifold in question
is endowed with a $3$-quasi-Sasakian structure. Although  the
notion of $3$-quasi-Sasakian structure is well known in literature,
it seems that a systematic study of these manifolds
has not been conducted so far. This is what we plan to do in the
present paper.

The first step in our work is to prove that the distribution
spanned by the three characteristic vector fields $\xi_1$,
$\xi_2$, $\xi_3$ of a $3$-quasi-Sasakian manifold is involutive
and defines a $3$-dimensional Riemannian and totally geodesic
foliation $\cal V$ of $M$. This property, in general, does not
hold for an almost $3$-contact metric manifold as we show by a
pair of simple examples. Taking into account
the geometry of the foliation $\cal V$
we show that $3$-quasi-Sasakian manifolds
divide into two classes: those manifolds for
which the foliation $\cal V$ has the local structure of an abelian
Lie group, and those for which $\cal V$ has the local structure of
the Lie group $SO(3)$ (or $SU(2)$).

For a $3$-quasi-Sasakian manifold one can consider the ranks of
the three  structures $(\phi_\alpha,\xi_\alpha,\eta_\alpha,g)$.
We will prove that these ranks coincide, allowing us
to  classify $3$-quasi-Sasakian manifolds according to
their well-defined rank, which is of the form $4p+1$ in
the abelian case and $4p+3$ in the non-abelian one.
Note that $3$-cosymplectic manifolds (rank $1$) and $3$-Sasakian manifolds
(rank $4n+3=\dim(M)$) are two representatives of each of the above
classes. Nevertheless we show examples of $3$-quasi-Sasakian
manifolds which are neither $3$-cosymplectic nor $3$-Sasakian. In
this sense we can say that the notion of $3$-quasi-Sasakian
manifold is a natural generalization of that of $3$-Sasakian and
that of $3$-cosymplectic manifold.

Furthermore, we prove a splitting theorem for any $3$-quasi-Sasakian
manifold $M$ assuming the integrability of one of the almost product
structures. We prove that if $M$ belongs to the class of
$3$-quasi-Sasakian manifolds with
$\left[\xi_\alpha,\xi_\beta\right]=2\xi_\gamma$, then $M$ is locally
the product of a $3$-Sasakian and a hyper-K\"{a}hlerian manifold,
whereas if $M$ belongs to the class of $3$-quasi-Sasakian manifolds
with $\left[\xi_\alpha,\xi_\beta\right]=0$, then $M$ is
$3$-cosymplectic. Finally, we find an application of the
integrability of the vertical distribution showing that $\cal V$ is
a minimum of the corrected energy in the sense of
Chac\'{o}n-Naveira-Weston \cite{chacon,chacon2} and Blair-Turgut
Vanli \cite{blair2, blair3}.

\section{Preliminaries}
By an \emph{almost contact manifold} we mean an odd-dimensional
manifold $M$ which carries a field $\phi$ of endomorphisms of the
tangent spaces, a vector field $\xi$, called \emph{characteristic}
or \emph{Reeb vector field}, and a $1$-form $\eta$ satisfying
$\phi^2=-I+\eta\otimes\xi$ and  $\eta\left(\xi\right)=1$, where $I
\colon TM\rightarrow TM$ is the identity mapping. From the
definition it follows also that $\phi\xi=0$, $\eta\circ\phi=0$ and
that the $(1,1)$-tensor field $\phi$ has constant rank $2n$ where $2n+1$
is the dimension of $M$ (cf. \cite{blair1}). Given an almost contact manifold
$(M,\phi,\xi,\eta)$ one can define an almost
complex structure $J$ on the product $M\times\mathbb{R}$ by setting
$J\left(X,f\frac{d}{dt}\right)=\left(\phi
X-f\xi,\eta\left(X\right)\frac{d}{dt}\right)$ for any
$X\in\Gamma\left(TM\right)$ and $f\in
C^{\infty}\left(M\times\mathbb{R}\right)$. Then the almost contact
manifold is said to be \emph{normal} if the almost
complex structure $J$ is integrable. The computation of the Nijenhuis
tensor of $J$ gives rise to the four tensors defined by
\begin{align*}
&N^{(1)}\left(X,Y\right)=\left[\phi,\phi\right]\left(X,Y\right)+2d\eta\left(X,Y\right)\xi,\\
&N^{(2)}\left(X,Y\right)=\left({\cal L}_{\phi
X}\eta\right)\left(Y\right)-\left({\cal L}_{\phi
Y}\eta\right)\left(X\right),\\
&N^{(3)}\left(X\right)=\left({\cal L}_{\xi}\phi\right)X,\\
&N^{(4)}\left(X\right)=\left({\cal L}_{\xi}\eta\right)\left(X\right),
\end{align*}
where $\left[\phi,\phi\right]$ is the Nijenhuis
tensor of $\phi$ and ${\cal L}_X$ denotes the Lie derivative with respect to the vector field $X$. One
finds that the structure $\left(\phi,\xi,\eta\right)$ is normal if
and only if $N^{(1)}$ vanishes identically; in particular,
$N^{(1)}=0$ implies the vanishing of $N^{(2)}$, $N^{(3)}$ and $N^{(4)}$ (cf.
\cite{sasaki2}).

Any almost contact manifold $\left(M,\phi,\xi,\eta\right)$ admits
a \emph{compatible} Riemannian metric $g$ such that $g\left(\phi
X,\phi Y\right) =
g\left(X,Y\right)-\eta\left(X\right)\eta\left(Y\right)$ for all
$X,Y\in\Gamma\left(TM\right)$. The manifold $M$ is then said to be
an \emph{almost contact metric manifold} with structure
$\left(\phi,\xi,\eta,g\right)$. The $2$-form $\Phi$ on $M$ defined
by $\Phi\left(X,Y\right)=g\left(X,\phi Y\right)$ is called the
\emph{fundamental $2$-form} of the almost contact metric manifold
$\left(M,\phi,\xi,\eta,g\right)$.

    Almost contact metric manifolds such that both $\eta$ and $\Phi$ are
closed are called \emph{almost cosymplectic manifolds} and those such that $d\eta=\Phi$ are called
\emph{contact metric manifolds}. Finally, a \emph{cosymplectic manifold} is a normal almost
cosymplectic manifold while a \emph{Sasakian manifold} is a normal contact metric manifold.

The notion of  quasi-Sasakian structure, introduced by D. E. Blair
in \cite{blair0}, unifies those of Sasakian and cosymplectic
structures. A \emph{quasi-Sasakian manifold} is defined as a
normal almost contact metric manifold whose fundamental $2$-form
is closed. A quasi Sasakian manifold $M$ is said to be of rank $2p$ (for some $p\leq n$) if
$\left(d\eta\right)^p\neq 0$ and $\eta\wedge\left(d\eta\right)^p=0$ on $M$,
and to be of rank $2p+1$ if $\eta\wedge\left(d\eta\right)^p\neq 0$ and
$\left(d\eta\right)^{p+1}=0$  on $M$ (cf. \cite{blair0,tanno}).
Blair proved that there are no quasi-Sasakian manifolds of even rank.
If the rank of $M$ is $2p+1$, then the module $\Gamma(TM)$ of vector
fields over $M$ splits into two submodules as follows:
$\Gamma(TM)={\cal E}^{2p+1}\oplus{\cal E}^{2q}$, $p+q=n$, where
\begin{equation*}
{\cal E}^{2q}=\{X\in\Gamma(TM)\; | \; i_X d\eta=0 \mbox{ and } i_X \eta=0\}
\end{equation*}
and  ${\cal E}^{2p+1}={\cal
E}^{2p}\oplus\left\langle\xi\right\rangle$, ${\cal E}^{2p}$ being
the orthogonal complement of ${\cal
E}^{2q}\oplus\left\langle\xi\right\rangle$ in
$\Gamma\left(TM\right)$. These modules satisfy $\phi {\cal
E}^{2p}={\cal E}^{2p}$ and $\phi {\cal E}^{2q}={\cal E}^{2q}$ (cf.
\cite{tanno}).

We will now mention some properties of quasi-Sasakian manifolds which
will be useful in the sequel.
\begin{lemma}[\cite{kanemaki1}]\label{teoremakanemaki}
A necessary and sufficient condition for an almost contact metric
manifold $\left(M,\phi,\xi,\eta,g\right)$ to be quasi-Sasakian is
that there exists a symmetric linear transformation field $A$ on
$M$, commuting with $\phi$, such that
\begin{equation}\label{formulakanemaki}
(\nabla_{X}\phi)Y=\eta(Y)AX-g(AX,Y)\xi,
\end{equation}
for any vector fields $X$ and $Y$ on $M$, where $\nabla$ denotes the Levi-Civita connection.
\end{lemma}

Replacing $Y$ with $\xi$ in \eqref{formulakanemaki} and applying $\phi$ we easily get
\begin{equation}\label{conseguenza1kanemaki}
\phi A = \nabla\xi,
\end{equation}
from which it follows that $\nabla_\xi\xi=\phi A\xi=A\phi\xi=0$, so
the integral curves of $\xi$ are geodesics. Moreover,  the
vanishing of $N^{(2)}$ implies
\begin{equation}
 d\eta\left(\phi X,Y\right)=d\eta\left(\phi Y,X\right)
    \label{phisimmetria}
\end{equation}
for all $X,Y\in\Gamma\left(TM\right)$. In quasi-Sasakian manifolds
the Reeb vector field has the following properties.

\begin{lemma}[\cite{blair0},\cite{olszak2}]
Let $\left(M,\phi,\xi,\eta,g\right)$ be a quasi-Sasakian manifold.
Then
\begin{description}
  \item[(i)] the Reeb vector field $\xi$ is Killing;
  \item[(ii)] the Ricci curvature in the direction of $\xi$ is given by
\begin{equation}\label{olszakricci}
\emph{Ric}\left(\xi\right)=\|\nabla\xi\|^2.
\end{equation}
\end{description}
\end{lemma}

We now come to the main topic of our paper, i.e. $3$-quasi-Sasakian geometry,
which is framed into the more general setting of almost $3$-contact geometry.
An  \emph{almost $3$-contact manifold}  is a
$\left(4n+3\right)$-dimensional smooth  manifold $M$ endowed with
three almost contact structures $\left(\phi_1,\xi_1,\eta_1\right)$,
$\left(\phi_2,\xi_2,\eta_2\right)$,
$\left(\phi_3,\xi_3,\eta_3\right)$ satisfying the following
relations, for any even permutation
$\left(\alpha,\beta,\gamma\right)$ of $\left\{1,2,3\right\}$,
\begin{gather}
\phi_\gamma=\phi_{\alpha}\phi_{\beta}-\eta_{\beta}\otimes\xi_{\alpha}=-\phi_{\beta}\phi_{\alpha}+\eta_{\alpha}\otimes\xi_{\beta},\\
\nonumber
\xi_{\gamma}=\phi_{\alpha}\xi_{\beta}=-\phi_{\beta}\xi_{\alpha}, \ \
\eta_{\gamma}=\eta_{\alpha}\circ\phi_{\beta}=-\eta_{\beta}\circ\phi_{\alpha}.
 \label{3-sasaki}
\end{gather}
This notion was introduced by Y. Y. Kuo \cite{kuo} and, independently,
by C. Udriste \cite{udriste}.  Kuo proved that given an
almost contact $3$-structure
$\left(\phi_\alpha,\xi_\alpha,\eta_\alpha\right)$, there exists a
Riemannian metric $g$ compatible with each of them and hence we can
speak of \emph{almost contact metric $3$-structures}. It is well
known that in any almost $3$-contact metric manifold the Reeb vector
fields $\xi_1,\xi_2,\xi_3$ are orthonormal with respect to the
compatible metric $g$ and that the structural group of the tangent
bundle is reducible to $Sp\left(n\right)\times I_3$. Moreover, by
putting
${\cal{H}}=\bigcap_{\alpha=1}^{3}\ker\left(\eta_\alpha\right)$ one
obtains a $4n$-dimensional distribution on $M$ and the tangent
bundle splits as the orthogonal sum $TM={\cal{H}}\oplus{\cal{V}}$,
where ${\cal V}=\left\langle\xi_1,\xi_2,\xi_3\right\rangle$. We will
call any vector belonging to the distribution $\cal H$ \emph{horizontal}
and any vector belonging to the distribution $\cal V$ \emph{vertical}.
An almost $3$-contact manifold $M$ is  said  to  be \emph{hyper-normal}
if each  almost contact structure $\left(\phi_\alpha,\xi_\alpha,\eta_\alpha\right)$ is
normal. It should be remarked that,  as it was proved in
\cite{yano1}, if two of the almost contact structures are normal,
then so is the third.

\begin{definition}
A \emph{$3$-quasi-Sasakian manifold} is a hyper-normal almost
$3$-contact metric manifold
$(M,\phi_\alpha,\xi_\alpha,\eta_\alpha,g)$ such that each
fundamental $2$-form $\Phi_\alpha$ is closed.
\end{definition}

The class of $3$-quasi-Sasakian manifolds includes as special cases the well-known
$3$-Sasakian and $3$-cosymplectic manifolds.

All manifolds considered here are assumed to be smooth i.e. of the class $\cal C^\infty$,
and connected; we denote by $\Gamma(\,\cdot\,)$ the set of all sections of a corresponding
bundle. We use the convention that $2u\wedge v=u\otimes v-v\otimes u$.

\section{The canonical foliation of a $3$-quasi-Sasakian manifold}

In this section we deal with the integrability of the distribution
$\cal V$ generated by the Reeb vector fields of an almost
$3$-contact manifold. We exhibit two examples of non-hyper-normal
almost $3$-contact structures. In the first one $\cal V$ is not
involutive, while in the second example we have an integrable
distribution.

\begin{example}
\emph{Let $\mathfrak g$ be the $7$-dimensional Lie algebra with basis
$\{X_1,X_2,X_3,X_4,$ $\xi_1,\xi_2,\xi_3\}$ and with the
Lie brackets defined by
$\left[X_i,X_j\right]=\left[X_i,\xi_\alpha\right]=0$ and
$\left[\xi_1,\xi_2\right]=\left[\xi_2,\xi_3\right]=\left[\xi_3,\xi_1\right]=X_1$,
for all $i,j\in\left\{1,2,3,4\right\}$ and
$\alpha\in\left\{1,2,3\right\}$. Let $G$ be a Lie group whose Lie
algebra is  $\mathfrak g$.  We define a left-invariant almost $3$-contact
structure $(\phi_\alpha,\xi_\alpha,\eta_\alpha)$ on $G$ putting $\phi_\alpha\xi_\beta=\epsilon_{\alpha\beta\gamma}\xi_\gamma$,
where $\epsilon_{\alpha\beta\gamma}$  is the totally antisymmetric symbol, and}
\begin{gather}\label{definizionephi7}\nonumber
\nonumber
\phi_1 X_1=X_2, \ \phi_1 X_2=-X_1, \ \phi_1 X_3=X_4, \ \phi_1 X_4=-X_3,\\
\phi_2 X_1=X_3, \ \phi_2 X_2=-X_4, \ \phi_2 X_3=-X_1, \ \phi_2
X_4=X_2,\\ \nonumber
 \phi_3 X_1=X_4, \ \phi_3 X_2=X_3, \ \phi_3
X_3=-X_2, \ \phi_3 X_4=-X_1,\nonumber
\end{gather}
\emph{and setting $\eta_\alpha\left(X_i\right)=0$, $\eta_\alpha\left(\xi_\beta\right)=\delta_{\alpha\beta}$, for all $\alpha,\beta\in\left\{1,2,3\right\}$ and $i\in\left\{1,2,3,4\right\}$. We note that this structure is not hyper-normal. Indeed, a straightforward computation yields $N_1^{(1)}\left(\xi_1,\xi_2\right)=-X_1+X_2\neq 0$ since $X_1$ and $X_2$ are linearly independent}.
\end{example}

\begin{example}
\emph{Using the above notations, consider the 7-dimensional Lie algebra  $\mathfrak
g'$ with brackets defined by
$[X_i,X_j]=[\xi_\alpha,\xi_\beta]=0$ and $[\xi_\alpha,
X_i]=\xi_\alpha$, for each $\alpha\in\left\{1,2,3\right\}$ and
$i,j\in\left\{1,2,3,4\right\}$. Let $G'$ be a Lie group whose Lie
algebra is $\mathfrak g'$ and let $\eta_\alpha$ and $\phi_\alpha$ be
the same tensors as above. Then
$(G',\phi_\alpha,\xi_\alpha,\eta_\alpha)$ is an almost $3$-contact
manifold which is not hyper-normal because
$N_\alpha^{(3)}(X_i)=[\xi_\alpha,\phi_\alpha
X_i]-\phi_\alpha[\xi_\alpha,X_i]=\pm\xi_\alpha$,
from which it follows that $N_\alpha^{(1)}\neq 0$.}
\end{example}

Therefore it makes sense to ask whether for a certain class of almost $3$-contact manifolds the distribution $\cal V$ is integrable. We will prove that this happens for a $3$-quasi-Sasakian manifold.
To this end we need the following Lemma that is a straightforward generalization of Lemma \ref{teoremakanemaki}.

\begin{lemma}\label{kanemaki3}
An almost $3$-contact metric manifold
$(M,\phi_\alpha,\xi_\alpha,\eta_\alpha,g)$ is a $3$-quasi-Sasakian
manifold if and only if it carries three symmetric tensor fields
$A_1$, $A_2$, $A_3$ such that, for each
$\alpha\in\left\{1,2,3\right\}$,
\begin{equation*}
(\nabla_{X}\phi_\alpha)Y=\eta_\alpha(Y)A_\alpha X-g(A_\alpha
X,Y)\xi_\alpha, \quad \quad \phi_\alpha A_\alpha=A_\alpha \phi_\alpha,
\end{equation*}
for all $X,Y\in\Gamma\left(TM\right)$.
\end{lemma}

\begin{theorem}\label{principale}
Let $(M,\phi_\alpha,\xi_\alpha,\eta_\alpha,g)$ be a
$3$-quasi-Sasakian manifold. Then the $3$-dimensional distribution
$\cal V$ generated by $\xi_1$, $\xi_2$, $\xi_3$ is integrable.
Moreover, $\cal V$ defines a totally geodesic and Riemannian
foliation of $M$.
\end{theorem}
\begin{proof}
First note that for any $\alpha\in\left\{1,2,3\right\}$ we have $\nabla_{\xi_\alpha}\xi_\alpha=0$,
each structure $(\phi_\alpha,\xi_\alpha,\eta_\alpha,g)$ being quasi-Sasakian.
Now, applying the first formula in Lemma \ref{kanemaki3} to $X=Y=\xi_\beta$  we obtain,
for any $\alpha \neq \beta$,
\begin{align*}
\nabla_{\xi_\beta} \left(\phi_\alpha \xi_\beta\right) - \phi_\alpha\left(\nabla_{\xi_\beta} \xi_\beta\right)=
\eta_\alpha \left(\xi_\beta\right) A_\alpha \xi_\beta- g\left(A_\alpha \xi_\beta,\xi_\beta\right)\xi_\alpha
\end{align*}
and hence
\begin{equation}
\nabla_{\xi_\beta}\xi_\gamma=-g\left(A_\alpha \xi_\beta,\xi_\beta\right)\xi_\alpha.
    \label{formulaintegrability}
\end{equation}
Thus, the distribution $\cal V$ is integrable with totally geodesic leaves.
Finally, $\cal V$ is a Riemannian foliation since the Reeb vector fields are Killing.
\end{proof}

\begin{remark}
\emph{Theorem \ref{principale} can be regarded as an analogue of a
result of Kuo and Tachibana \cite{kuo2} stating that the
distribution spanned by the three Reeb vector fields of a
hyper-normal $3$-contact manifold is integrable with totally
geodesic leaves.}
\end{remark}

\begin{corollary}\label{preserva}
In any $3$-quasi-Sasakian manifold the operators $A_\alpha$ preserve the vertical and horizontal distributions.
\end{corollary}
\begin{proof}
Applying the first formula in Lemma \ref{kanemaki3} to $X=Y=\xi_\alpha$ one easily gets $A_\alpha \xi_\alpha = g\left(A_\alpha \xi_\alpha,\xi_\alpha\right)\xi_\alpha$. For $\alpha\neq\beta$, applying the formula $A_\alpha \phi_\alpha=\nabla \xi_\alpha$ and \eqref{formulaintegrability} we obtain
\begin{equation*}
    A_\alpha \xi_\beta = -A_\alpha \phi_\alpha \xi_\gamma =-\nabla_{\xi_\gamma} \xi_\alpha = g\left(A_\beta \xi_\gamma ,\xi_\gamma\right) \xi_\beta.
\end{equation*}
Thus, the distribution $\cal V$ is preserved by $A_\alpha$. Moreover, since the operators $A_\alpha$ are symmetric  they also preserve the horizontal distribution $\cal H$.
\end{proof}

\begin{corollary}\label{facile2}
Let $(M,\phi_\alpha,\xi_\alpha,\eta_\alpha,g)$ be a
$3$-quasi-Sasakian manifold. Then,  for all $X\in\Gamma\left(\cal
H\right)$ and for all $\alpha,\beta\in\left\{1,2,3\right\}$,
$d\eta_\alpha\left(X,\xi_\beta\right)=0$.
\end{corollary}
\begin{proof}
For $\alpha=\beta$ the result is proved in \cite{blair0}. Next, let $\alpha\neq\beta$.
Using the integrability of the distribution
${\cal V}$ and the fact that $\xi_\beta$ is  Killing  we have
$2d\eta_{\alpha}\left(X,\xi_\beta\right)=-\eta_\alpha\left(\left[X,\xi_\beta\right]\right)=-g\left(\left[X,\xi_\beta\right],\xi_\alpha\right)=\xi_\beta\left(g\left(X,\xi_\alpha\right)\right)-g\left(X,\left[\xi_\beta,\xi_\alpha\right]\right)=0$.
\end{proof}

\begin{corollary}\label{facile3}
In any $3$-quasi-Sasakian manifold each Reeb vector field
$\xi_\alpha$ is an infinitesimal automorphism with respect to the
horizontal distribution $\cal H$.
\end{corollary}
\begin{proof}
From Corollary \ref{facile2} it follows that for any $X\in \Gamma\left(\cal H\right)$ and any $\alpha,\beta\in\{1,2,3\}$,
$0=d\eta_\alpha\left(X,\xi_\beta\right)=-\frac{1}{2}\eta_\alpha\left([X,\xi_\beta]\right)$ and hence $[X,\xi_\beta]\in \Gamma\left(\cal H\right)$.
\end{proof}

\section{Lie group structures on the leaves of ${\cal V}$}

In this section we prove that the canonical foliation $\cal V$ of a
$3$-quasi-Sasakian manifold has locally either the structure of the orthogonal group
or that of an abelian group. Accordingly, $3$-quasi-Sasakian manifolds split into two main classes.

\begin{lemma}\label{giulia}
Let $(M,\phi_\alpha,\xi_\alpha,\eta_\alpha,g)$ be a hyper-normal
almost 3-contact metric manifold. Then, for any even permutation
$\left(\alpha,\beta,\gamma\right)$ of $\left\{1,2,3\right\}$ we have
\begin{equation*}
[\xi_\alpha,\xi_\beta]=[\xi_\alpha,\xi_\beta]_{\cal H}+f\xi_\gamma,
\end{equation*}
where $f$ is the function given by
$f=\eta_\gamma([\xi_\alpha,\xi_\beta])=\eta_\alpha([\xi_\beta,\xi_\gamma])=\eta_\beta([\xi_\gamma,\xi_\alpha])$.
\end{lemma}
\begin{proof}
Since $M$ is hyper-normal we have that ${\cal L}_{\xi_\alpha}\phi_\alpha=0$.
The computation of $({\cal L}_{\xi_\alpha}\phi_\alpha)\xi_\gamma=0$ yields
$[\xi_\alpha,\xi_\beta]=-\phi_\alpha[\xi_\alpha,\xi_\gamma]$,
while $({\cal L}_{\xi_\beta}\phi_\beta)\xi_\gamma=0$ yields
$[\xi_\beta,\xi_\alpha]=\phi_\beta[\xi_\beta,\xi_\gamma]$.
It follows that
$[\xi_\alpha,\xi_\beta]\in\ker(\eta_\alpha)\cap\ker(\eta_\beta)$.
Thus
\begin{equation}
[\xi_\alpha,\xi_\beta]=[\xi_\alpha,\xi_\beta]_{\cal H}+f_\gamma\xi_\gamma,\label{partenza}
\end{equation}
for some functions $f_\gamma$. We prove that $f_1=f_2=f_3$.
Applying $\phi_\alpha$ to  \eqref{partenza} we get
\begin{equation*}
\phi_\alpha[\xi_\alpha,\xi_\beta]=\phi_\alpha[\xi_\alpha,\xi_\beta]_{\cal H}-f_\gamma\xi_\beta.
\end{equation*}
On the other hand, from $({\cal L}_{\xi_\alpha}\phi_\alpha)\xi_\beta=0$ it follows that
\begin{equation*}
\phi_\alpha[\xi_\alpha,\xi_\beta]=[\xi_\alpha,\xi_\gamma]=-[\xi_\gamma,\xi_\alpha]_{\cal
H}-f_\beta\xi_\beta.
\end{equation*}
Hence, $\phi_\alpha[\xi_\alpha,\xi_\beta]_{\cal
H}-f_\gamma\xi_\beta=-[\xi_\gamma,\xi_\alpha]_{\cal H}-f_\beta\xi_\beta$, from which it follows
that $f_\beta=f_\gamma$.
\end{proof}
\\

Now we prove that for a $3$-quasi-Sasakian manifold the function $f$
which appears in Lemma \ref{giulia} is necessarily constant.

\begin{theorem}\label{classification}
Let $(M,\phi_\alpha,\xi_\alpha,\eta_\alpha,g)$ be a
$3$-quasi-Sasakian manifold. Then, for any even permutation
$(\alpha,\beta,\gamma)$ of $\left\{1,2,3\right\}$ and for some $c\in \mathbb R$
\begin{equation*}
\left[\xi_\alpha,\xi_\beta\right]=c\xi_\gamma.
\end{equation*}
\end{theorem}
\begin{proof}
From Lemma  \ref{giulia} and Theorem \ref{principale} it follows
that $\left[\xi_\alpha,\xi_\beta\right]=f\xi_\gamma$, where $f\in
C^{\infty}\left(M\right)$ is the function given by
$f=\eta_\alpha\left(\left[\xi_\beta,\xi_\gamma\right]\right)=-2d\eta_\alpha\left(\xi_\beta,\xi_\gamma\right)$.
We have thus only to prove that such a function is constant. Indeed,
the vanishing of $N^{(4)}$ implies that ${\cal
L}_{\xi_\alpha}d\eta_\alpha=0$. Then, for any even permutation
$(\alpha,\beta,\gamma)$ of $\left\{1,2,3\right\}$,
\begin{align*}
0&=\left({\cal
L}_{\xi_\alpha}d\eta_\alpha\right)\left(\xi_\beta,\xi_\gamma\right)\\
&=\xi_\alpha\left(d\eta_\alpha\left(\xi_\beta,\xi_\gamma\right)\right)-d\eta_\alpha\left(\left[\xi_\alpha,\xi_\beta\right],\xi_\gamma\right)-d\eta_\alpha\left(\xi_\beta,\left[\xi_\alpha,\xi_\gamma\right]\right)\\
&=-\frac{1}{2}\xi_\alpha\left(f\right)-fd\eta_\alpha\left(\xi_\gamma,\xi_\gamma\right)+fd\eta_\alpha\left(\xi_\beta,\xi_\beta\right)\\
&=-\frac{1}{2}\xi_\alpha\left(f\right),
\end{align*}
so that $f$ is constant along the leaves of $\cal V$. It remains to
show that $X\left(f\right)=0$ for all $X\in\Gamma\left(\cal
H\right)$. In fact, using the formula of the differential of a
$2$-form, we find
\begin{align*}
-\frac{1}{2}X\left(f\right)&=X\left(d\eta_\alpha\left(\xi_\beta,\xi_\gamma\right)\right)\\
&=3d^2\eta_\alpha\left(X,\xi_\beta,\xi_\gamma\right)-\xi_\beta\left(d\eta_\alpha\left(\xi_\gamma,X\right)\right)-\xi_\gamma\left(d\eta_\alpha\left(X,\xi_\beta\right)\right)\\
&\quad+d\eta_\alpha\left(\left[X,\xi_\beta\right],\xi_\gamma\right)
+d\eta_\alpha\left(\left[\xi_\beta,\xi_\gamma\right],X\right)+d\eta_\alpha\left(\left[\xi_\gamma,X\right],\xi_\beta\right)\\
&=fd\eta_\alpha\left(\xi_\alpha,X\right)=0,
\end{align*}
where we have used Corollary \ref{facile2}, Corollary \ref{facile3}
and the integrability of $\cal V$.
\end{proof}
\\

Using Theorem \ref{classification} we may therefore divide
$3$-quasi-Sasakian manifolds in two classes according
to the behaviour of the leaves of the canonical
foliation $\cal V$: those $3$-quasi-Sasakian manifolds for which
each leaf of $\cal V$ is locally $SO\left(3\right)$ (or
$SU\left(2\right)$) (which corresponds to take in Theorem
\ref{classification} the constant $c\neq 0$), and those for which
each leaf of $\cal V$ is locally an abelian group (this corresponds
to the case $c=0$).

Note that $3$-Sasakian manifolds and $3$-cosymplectic manifolds are
representatives of each of the above classes. However we emphasize
the circumstance that $3$-Sasakian and $3$-cosymplectic manifolds do
not exhaust the above two classes. This can be seen in the
construction of the following example.

\begin{example}
\emph{Let us denote  the canonical global coordinates on $\mathbb{R}^{4n+3}$ by \hfil
 $x_1,\ldots,x_n$, $y_1,\ldots,y_n, u_1,\ldots,u_n,v_1,\ldots,v_n,z_1,z_2,z_3$. We consider the open submanifold $M$ of
$\mathbb{R}^{4n+3}$ obtained by removing the points where
$\sin\left(z_2\right)=0$ and define three vector fields and three
$1$-forms on $M$ by
\begin{align*}
\xi_1&=c\frac{\partial}{\partial z_1}, \\
\xi_2&=c\left(\cos\left(z_1\right)\cot\left(z_2\right)\frac{\partial}{\partial
z_1}+\sin\left(z_1\right)\frac{\partial}{\partial
z_2}-\frac{\cos\left(z_1\right)}{\sin\left(z_2\right)}\frac{\partial}{\partial
z_3}\right), \\
\xi_3&=c\left(-\sin\left(z_1\right)\cot\left(z_2\right)\frac{\partial}{\partial
z_1}+\cos\left(z_1\right)\frac{\partial}{\partial
z_2}+\frac{\sin\left(z_1\right)}{\sin\left(z_2\right)}\frac{\partial}{\partial
z_3}\right)
\end{align*}
for some non-zero real number $c$, and
\begin{align*}
\eta_1&=\frac{1}{c}\left(dz_1+\cos\left(z_2\right)dz_3\right),\\
\eta_2&=\frac{1}{c}\left(\sin\left(z_1\right)dz_2-\cos\left(z_1\right)\sin\left(z_2\right)dz_3\right),\\
\eta_3&=\frac{1}{c}\left(\cos\left(z_1\right)dz_2+\sin\left(z_1\right)\sin\left(z_2\right)dz_3\right).
\end{align*}
A simple computation shows that $\left[\xi_{\alpha},\xi_{\beta}\right]=c\xi_{\gamma}$ for any cyclic permutation, and $\eta_\alpha\left(\xi_\beta\right)=\delta_{\alpha\beta}$ for any $\alpha,\beta\in\left\{1,2,3\right\}$. Now, in order to define three
tensor fields $\phi_\alpha$ and a Riemannian metric $g$ on $M$, we
set, for each $i\in\left\{1,\ldots,n\right\}$,
$X_i=\frac{\partial}{\partial x_i}, \ Y_i=\frac{\partial}{\partial
y_i}, \ U_i=\frac{\partial}{\partial u_i}, \
V_i=\frac{\partial}{\partial v_i}$.
Let $g$ be the Riemannian metric with respect to which
$\{X_1,\ldots,X_n,Y_1,\ldots,Y_n,$ \ $U_1,\ldots,U_n$,
$V_1,\ldots,V_n,\xi_1,\xi_2,\xi_3\}$ is a (global) orthonormal
frame, and $\phi_1$, $\phi_2$, $\phi_3$ be the tensor fields defined
by putting $\phi_\alpha\xi_\beta=\epsilon_{\alpha\beta\gamma}\xi_\gamma$ and, for all $i\in\left\{1,\ldots,n\right\}$,
\begin{gather}\label{definizionephi}\nonumber
\nonumber
\phi_1 X_i=Y_i, \ \phi_1 Y_i=-X_i, \ \phi_1 U_i=V_i, \ \phi_1 V_i=-U_i,\\
\phi_2 X_i=U_i, \ \phi_2 Y_i=-V_i, \ \phi_2 U_i=-X_i, \ \phi_2
V_i=Y_i,\\ \nonumber
 \phi_3 X_i=V_i, \ \phi_3 Y_i=U_i, \ \phi_3
U_i=-Y_i, \ \phi_3 V_i=-X_i.\nonumber
\end{gather}
Note that we have also $\eta_\alpha\left(X_i\right) =\eta_\alpha\left(Y_i\right) =\eta_\alpha\left(U_i\right) =\eta_\alpha\left(V_i\right)=0$.
By a lengthy computation that we omit one can verify that $\left(\phi_\alpha,\xi_\alpha,\eta_\alpha,g\right)$
is an almost contact metric $3$-structure.  Observe that the structure is not $3$-cosymplectic
since the Reeb vector fields do not commute. Nevertheless it is neither $3$-Sasakian, since it admits a Darboux-like coordinate system (cf. \cite{cappellettidenicola}). Some simple computations show that  $\left(\phi_\alpha,\xi_\alpha,\eta_\alpha,g\right)$ is hyper-normal. For instance,
\begin{align*}
N^{(1)}_{1}\left(\xi_{2},\xi_{3}\right)&= \phi_{1}^2\left[\xi_{2},\xi_{3}\right] + \left[\phi_{1}\xi_{2},\phi_{1}\xi_{3} \right]-\phi_{1}\left[\xi_{2},\phi_{1}\xi_{3} \right] - \phi_{1}\left[\phi_{1}\xi_{2},\xi_{3} \right] - \eta_1\left(\left[\xi_{2},\xi_{3} \right]\right)\xi_{1}\\
    &= c\phi_{1}^2\xi_{1} - \left[\xi_3,\xi_2\right] + \phi_{1}\left[\xi_2 ,\xi_2 \right] - \phi_{1}\left[\xi_3 ,\xi_3 \right] - c\eta_1\left(\xi_{1}\right)\xi_{1}=0.
\end{align*}
Furthermore, the fundamental $2$-forms $\Phi_\alpha$ are closed. Indeed,  for instance we have
\begin{align*}
3d\Phi_1\left(X_i,X_j,X_k\right) &= X_i(\Phi_1\left(X_j,X_k\right)) - X_j(\Phi_1\left(X_i,X_k\right)) +X_k(\Phi_1\left(X_i,X_j\right))\\
&= X_i (g\left(X_j,\phi_1 X_k\right)) - X_j (g\left(X_i,\phi_1 X_k\right)) + X_k (g\left(X_i,\phi_1 X_k\right))=0,
\end{align*}
since the scalar products in the above expression are constant functions. Thus,
$\left(M,\phi_\alpha,\xi_\alpha,\eta_\alpha,g\right)$ is a
$3$-quasi-Sasakian manifold. Furthermore,  one can
show that $M$ is $\eta$-Einstein, its Ricci tensor being given by
$\textrm{Ric}=\frac{c^2}{2}\left(\eta_1\otimes\eta_1+\eta_2\otimes\eta_2+\eta_3\otimes\eta_3\right)$.
Thus, differently from $3$-Sasakian and $3$-cosymplectic geometry,
there are $3$-quasi-Sasakian manifolds which are not Einstein.}
\end{example}

\section{The rank of a $3$-quasi-Sasakian manifold}
For a $3$-quasi-Sasakian manifold one can consider the ranks of the three  structures $(\phi_\alpha,\xi_\alpha,\eta_\alpha,g)$.
In this section we will prove that these ranks coincide allowing us  to classify $3$-quasi-Sasakian manifolds. First we need some technical results.
\begin{lemma}
Let $(M,\phi_\alpha,\xi_\alpha,\eta_\alpha,g)$ be an almost $3$-contact metric manifold. Then we have
\begin{equation}\label{formulagenerale}
\eta_\alpha=\frac{1}{2}i_{\xi_\beta}\Phi_\gamma,
\end{equation}
for an even permutation $(\alpha,\beta,\gamma)$ of $\left\{1,2,3\right\}$.
Furthermore, if the fundamental $2$-forms $\Phi_\alpha$ are closed, we have
\begin{equation}\label{formuladifferenziale}
d\eta_\alpha=\frac{1}{2}{\cal L}_{\xi_\beta}\Phi_\gamma.
\end{equation}
\end{lemma}
\begin{proof}
For any $X\in\Gamma\left(TM\right)$ we have
$\eta_\alpha\left(X\right)=g\left(X,\xi_\alpha\right)=-g\left(X,\phi_\gamma\xi_\beta\right)=-\Phi_\gamma\left(X,\xi_\beta\right)=\frac{1}{2}i_{\xi_\beta}\Phi_\gamma\left(X\right)$.
Assuming $d\Phi_\alpha=0$ and applying the Cartan formula we get ${\cal
L}_{\xi_\beta}\Phi_\gamma=di_{\xi_\beta}
\Phi_\gamma+i_{\xi_\beta}d\Phi_\gamma=di_{\xi_\beta}\Phi_\gamma=2d\eta_\alpha$.
\end{proof}

\begin{lemma}\label{lemmarango}
In any $3$-quasi-Sasakian manifold
$(M,\phi_\alpha,\xi_\alpha,\eta_\alpha,g)$ we have, for any even permutation
$\left(\alpha,\beta,\gamma\right)$ of $\left\{1,2,3\right\}$ and
for all $X,Y\in\Gamma\left(\cal H\right)$,
    $$d\eta_\alpha\left(\phi_\beta X,Y\right)=d\eta_\gamma\left(X,Y\right).$$
\end{lemma}
\begin{proof}
First \ of \ all \ note \ that \ $\Phi_\gamma\left(\phi_\beta
X,Y\right)=g\left(\phi_\beta X,\phi_\gamma
Y\right)=-g\left(X,\phi_\beta\phi_\gamma Y\right)=$ \
$-g\left(X,\phi_\alpha Y\right)=-\Phi_\alpha\left(X,Y\right)$. Using
this and \eqref{formuladifferenziale} we get
\begin{align*}
2d\eta_\alpha\left(\phi_\beta X,Y\right)&=\left({\cal
L}_{\xi_\beta}\Phi_\gamma\right)\left(\phi_\beta X,Y\right)\\
&=\xi_\beta\left(\Phi_\gamma\left(\phi_\beta
X,Y\right)\right)-\Phi_\gamma\left(\left[\xi_\beta,\phi_\beta
X\right],Y\right)-\Phi_\gamma\left(\phi_\beta
X,\left[\xi_\beta,Y\right]\right)\\
&=\xi_\beta\left(\Phi_\gamma\left(\phi_\beta
X,Y\right)\right)-\Phi_\gamma\left(\phi_\beta\left[\xi_\beta,X\right],Y\right)-\Phi_\gamma\left(\phi_\beta
X,\left[\xi_\beta,Y\right]\right)\\
&=-\xi_\beta\left(\Phi_\alpha\left(X,Y\right)\right)+\Phi_\alpha\left(\left[\xi_\beta,X\right],Y\right)+\Phi_\alpha\left(X,\left[\xi_\beta,Y\right]\right)\\
&=-\left({\cal L}_{\xi_\beta}\Phi_\alpha\right)\left(X,Y\right)\\
&=2d\eta_\gamma\left(X,Y\right),
\end{align*}
where we have used the fact that $N^{(3)}_\beta=0$.
\end{proof}
\\

The next lemma is an immediate consequence of Lemma \ref{lemmarango} and Corollary \ref{facile2}.
\begin{lemma}\label{lemmarango2}
In any $3$-quasi-Sasakian manifold
$(M,\phi_\alpha,\xi_\alpha,\eta_\alpha,g)$ we have, for any $\alpha\neq\beta$
 and
for all $X,Y\in\Gamma\left(\cal H\right)$,
\begin{enumerate}
    \item[\emph{(}i\emph{)}] $d\eta_\alpha\left(X,\phi_\alpha Y\right)=d\eta_\beta\left(X,\phi_\beta Y\right)$,
    \item[\emph{(}ii\emph{)}] $d\eta_\alpha\left(\phi_\beta X,\phi_\beta Y\right)=-d\eta_\alpha\left(X,Y\right).$
\end{enumerate}
\end{lemma}

\begin{lemma}\label{ixetaalpha}
In any $3$-quasi-Sasakian manifold
$(M,\phi_\alpha,\xi_\alpha,\eta_\alpha,g)$ we have
$i_X d\eta_\alpha=0$ if and only if $i_X d\eta_\beta=0$ for any
$X\in\Gamma(\cal H)$ and $\alpha,\beta\in\{1,2,3\}$.
\end{lemma}
\begin{proof}
Assume $i_X d\eta_\alpha=0$. Then, from Lemma \ref{lemmarango} it follows that
$i_X d\eta_\beta(Y)=-2d\eta_\beta(Y,X) =2d\eta_\alpha(\phi_\gamma Y,X)=-i_X d\eta_\alpha(\phi_\gamma Y)=0$.
\end{proof}

\begin{theorem}\label{rango}
Let $(M,\phi_\alpha,\xi_\alpha,\eta_\alpha,g)$ be a
$3$-quasi-Sasakian manifold of dimension $4n+3$. Then the $1$-forms $\eta_1$, $\eta_2$
and $\eta_3$ have the same rank $4l+3$
or $4l+1$, for some $l\leq n$, according to
$\left[\xi_\alpha,\xi_\beta\right]=c\xi_\gamma$ with $c\neq 0$, or
$\left[\xi_\alpha,\xi_\beta\right]=0$, respectively.
\end{theorem}
\begin{proof}
Let us consider the quasi-Sasakian structures
$(\phi_\alpha,\xi_\alpha,\eta_\alpha,g)$ and assume that the 1-forms
$\eta_\alpha$ have ranks $2p_\alpha+1$. Then, according to \cite{tanno}, we have three
decompositions
\[\Gamma(TM)={\cal E}^{2p_\alpha+1}\oplus{\cal E}^{2q_\alpha},\]
where
\[{\cal E}^{2q_\alpha}=\{X\in\Gamma(TM)\; | \; i_{X}d\eta_\alpha=0 \mbox{ and } i_{X}\eta_\alpha=0\}.\]

Let us consider the class of 3-quasi-Sasakian manifolds such that
$[\xi_\alpha,\xi_\beta]=c\xi_\gamma$ with $c\neq 0$. In this case
we have $i_{\xi_\gamma}d\eta_\alpha(\xi_\beta)
=-i_{\xi_\beta}d\eta_\alpha(\xi_\gamma)=c$, which
implies that $\xi_\beta$, $\xi_\gamma\not\in{\cal E}^{2q_\alpha}$ and
\[{\cal E}^{2q_\alpha}=\{X\in\Gamma({\cal H})\; | \; i_X d\eta_\alpha=0\}.\]
Hence, by virtue of Lemma \ref{ixetaalpha}, ${\cal E}^{2q_\alpha}={\cal E}^{2q_\beta}$ and $q_\alpha=q_\beta$ for
$\alpha,\beta\in\{1,2,3\}$. It can be easily
seen that $\phi_\alpha {\cal E}^{2q}= {\cal E}^{2q}$ for $\alpha\in\{1,2,3\}$, where we have set $q=q_\alpha$.
Furthermore, the restrictions $\phi_\alpha|_{{\cal E}^{2q}}$ of the mappings $\phi_\alpha$ to ${\cal E}^{2q}$ define a quaternionic structure in the submodule ${\cal E}^{2q}$ which will be denoted by ${\cal E}^{4m}$ since its dimension is a multiple of 4.
It follows that in this case the rank of $M$ is $4l+3$, where $l=n-m$.

Let us now suppose $[\xi_\alpha,\xi_\beta]=0$. In this case
$\xi_\beta\in{\cal E}^{2q_\alpha}$ for $\beta\ne\alpha$. Indeed,
by Corollary \ref{facile2},
$i_{\xi_\beta}d\eta_\alpha(X)=0$ for any $X\in\Gamma({\cal H})$.
Moreover, $i_{\xi_\beta}d\eta_\alpha(\xi_\gamma)=0$ for any
$\gamma\in\{1,2,3\}$. It follows that
\[{\cal E}^{2q_\alpha}=\left\langle\xi_\beta,\xi_\gamma\right\rangle \oplus \{X\in\Gamma({\cal H})\; | \;
i_X d\eta_\alpha=0\}.\] This implies again that $q_\alpha=q_\beta$.
Moreover, the restrictions of the mappings $\phi_\alpha$ define as  before a quaternionic structure in
$\{X\in\Gamma({\cal H})\; | \; i_X d\eta_\alpha=0\}$, whose dimension is then $4m$ for some $m\leq n$. In this case the rank of the manifold
will be $4l+1$, where $l=n-m$.
\end{proof}
\\

According \ to \ Theorem \ref{rango}, \ we \ say \ that \ a \
$3$-quasi-Sasakian \
 manifold $(M,\phi_\alpha,\xi_\alpha,\eta_\alpha,g)$ has \emph{rank}
 $4l+3$ or $4l+1$ if any quasi-Sasakian structure has such rank. We may thus classify $3$-quasi-Sasakian manifolds of dimension
$4n+3$, according to their rank. For any $l\in\{0,\ldots,n\}$ we
have one class of manifolds such that
$\left[\xi_\alpha,\xi_\beta\right]=c\xi_\gamma$ with $c\neq 0$,
and one class of manifolds with
$\left[\xi_\alpha,\xi_\beta\right]=0$. The total number of classes
amounts then to $2n+2$.

In the following we will use the notation ${\cal E}^{4m}:=\{X\in\Gamma({\cal H})\; | \; i_X d\eta_\alpha=0\}$, while
${\cal E}^{4l}$ will be the orthogonal complement of ${\cal E}^{4m}$ in $\Gamma({\cal H})$,
${\cal E}^{4l+3}:={\cal E}^{4l} \oplus \Gamma({\cal V})$, and
${\cal E}^{4m+3}:={\cal E}^{4m} \oplus \Gamma({\cal V})$.

We now consider the class of $3$-quasi-Sasakian manifolds such that
$[\xi_\alpha,\xi_\beta]=c\xi_\gamma$ with $c\neq 0$ and let $4l+3$ be the rank.
In this case, according to \cite{blair0}, we define for each structure
$(\phi_\alpha,\xi_\alpha,\eta_\alpha,g)$ two $(1,1)$-tensor fields
$\psi_\alpha$ and $\theta_\alpha$ by putting
\begin{equation*}
 \psi_\alpha X=\left\{
             \begin{array}{ll}
               \phi_\alpha X, & \hbox{if $X\in {\cal{E}}^{4l+3}$;}\\
               0, & \hbox{if $X\in {\cal{E}}^{4m}$;}
             \end{array}
              \right.
\ \textrm{  } \
 \theta_\alpha X=\left\{
             \begin{array}{ll}
               0, & \hbox{if $X\in {\cal{E}}^{4l+3}$;}\\
               \phi_\alpha X, & \hbox{if $X\in {\cal{E}}^{4m}$.}\\
             \end{array}
              \right.
\end{equation*}
Note that, for each $\alpha\in\left\{1,2,3\right\}$ we have
$\phi_\alpha=\psi_\alpha+\theta_\alpha$.
Next, we define a new (pseudo-Riemannian, in general) metric
$\bar{g}$ on $M $ setting
\begin{equation*}
\bar{g}\left(X,Y\right)=\left\{
                          \begin{array}{ll}
                            -d\eta_\alpha\left(X,\phi_\alpha Y\right), & \hbox{for $X,Y\in{\cal E}^{4l}$;} \\
                            g\left(X,Y\right), & \hbox{elsewhere.}
                          \end{array}
                        \right.
\end{equation*}
This definition is well posed by virtue of \eqref{phisimmetria} and of Lemma
\ref{lemmarango2}. $(M,\phi_\alpha,\xi_\alpha,\eta_\alpha,\bar{g})$
is in fact a hyper-normal almost $3$-contact metric manifold, in general
non-$3$-quasi-Sasakian.
We are now able to prove the following decomposition theorem.
\begin{theorem}
Let $(M^{4n+3},\phi_\alpha,\xi_\alpha,\eta_\alpha,g)$ be a
$3$-quasi-Sasakian manifold of rank $4l+3$ with $\left[\xi_\alpha,\xi_\beta\right]=2\xi_\gamma$. Assume $[\theta_\alpha, \theta_\alpha]=0$ for some $\alpha\in\{1,2,3\}$ and $\bar{g}$ positive definite on ${\cal E}^{4l}$.
Then $M^{4n+3}$ is locally the  product of a $3$-Sasakian manifold $M^{4l+3}$ and a hyper-K\"{a}hlerian manifold $M^{4m}$ with $m=n-l$.
\end{theorem}
\begin{proof}
The almost product structure defined by $-\theta_\alpha^2$ and $-\psi_\alpha^2+ \eta_\alpha\otimes\xi_\alpha$ is integrable (see e.g. \cite{yano0}), since $[\theta_\alpha,\theta_\alpha]=0$ implies $[-\theta_\alpha^2,-\theta_\alpha^2]=0$. Thus, locally $M$ is the product of the manifolds $M^{4l+3}$ and $M^{4m}$ whose localized modules of vector fields are ${\cal E}^{4l+3}$ and ${\cal E}^{4m}$, respectively. Since $\psi_\alpha$ and $\phi_\alpha$ agree on ${\cal E}^{4l+3}$, $(\psi_\alpha,\xi_\alpha,\eta_\alpha)|_{{\cal E}^{4l+3}}$ is an almost contact $3$-structure. Furthermore, the metric $\bar{g}$ is by definition compatible with the three almost contact structures and $d\eta_\alpha=\bar\Phi_\alpha$ on ${\cal E}^{4l+3}$, so  $(\psi_\alpha,\xi_\alpha,\eta_\alpha,\bar{g})|_{{\cal E}^{4l+3}}$ is a contact metric $3$-structure over $M^{4l+3}$ and thus it is $3$-Sasakian, by a result of Kashiwada \cite{kashiwada}. Since $\theta_\alpha$ agrees with $\phi_\alpha$ on ${\cal E}^{4m}$, the maps $\theta_\alpha|_{{\cal E}^{4m}}$ define a quaternionic structure which is compatible with the metric $\bar{g}|_{{\cal E}^{4m}}$.  Finally, define the $2$-forms $\bar{\Theta}_\alpha$ by $\bar{\Theta}_\alpha\left(X,Y\right)=\bar{g}\left(X,\theta_\alpha Y\right)$ for any $X,Y\in {{\cal E}^{4m}}$. We have $\bar{\Theta}_\alpha=\bar{\Phi}_\alpha|_{{\cal E}^{4m}}$ and hence $d\bar{\Theta}_\alpha=0$.  By virtue of Hitchin Lemma \cite{hitchin} the structure defined on $M^{4m}$ turns out to be hyper-K\"{a}hlerian.
\end{proof}
\\

We now consider the class of $3$-quasi-Sasakian manifolds such that
$[\xi_\alpha,\xi_\beta]=0$ and let $4l+1$ be the rank.
In this case we define for each structure
$(\phi_\alpha,\xi_\alpha,\eta_\alpha,g)$ two $(1,1)$-tensor fields
$\psi_\alpha$ and $\theta_\alpha$ by putting
\begin{equation*}
 \psi_\alpha X=\left\{
             \begin{array}{ll}
               \phi_\alpha X, & \hbox{if $X\in {\cal{E}}^{4l}$;}\\
               0, & \hbox{if $X\in {\cal{E}}^{4m+3}$;}
             \end{array}
              \right.
\ \textrm{  } \
 \theta_\alpha X=\left\{
             \begin{array}{ll}
               0, & \hbox{if $X\in {\cal{E}}^{4l}$;}\\
               \phi_\alpha X, & \hbox{if $X\in {\cal{E}}^{4m+3}$.}\\
             \end{array}
              \right.
\end{equation*}

Note that for each $\alpha$ the maps $-\psi_\alpha^2$ and $-\theta_\alpha^2+ \eta_\alpha\otimes\xi_\alpha$ define an almost product structure which is integrable if and only if $[-\psi_\alpha^2,-\psi_\alpha^2]=0$ or, equivalently, $[\psi_\alpha,\psi_\alpha]=0$. Under this assumption the structure turns out to be $3$-cosymplectic as it is shown by the following theorem.

\begin{theorem}
Let \ $(M,\phi_\alpha,\xi_\alpha,\eta_\alpha,g)$  be  a
$3$-quasi-Sasakian  manifold of rank $4l+1$ such  that
$\left[\xi_\alpha,\xi_\beta\right]=0$ for any
$\alpha,\beta\in\{1,2,3\}$ and $[\psi_\alpha, \psi_\alpha]=0$ for
some $\alpha\in\{1,2,3\}$. Then $M$ is  a $3$-cosymplectic
manifold.
\end{theorem}
\begin{proof}
From $[\psi_\alpha, \psi_\alpha]=0$ it follows that the
corresponding quasi-Sasakian structure
$(\phi_\alpha,\xi_\alpha,\eta_\alpha,g)$ is cosymplectic, i.e.
$d\eta_\alpha=0$ (see \cite{blair0}, page 339). Then, for
any $\beta\neq \alpha$ and any $X\in\Gamma(\cal H)$, $Y\in\Gamma(TM)$
by Lemma \ref{ixetaalpha} we obtain $2d\eta_\beta\left(X,Y\right)=i_X d\eta_\beta\left(Y\right)=0$.
If $X=\xi_\delta$, $Y=\xi_{\rho}$ for some $\delta,\rho\in\{1,2,3\}$, then
$2d\eta_\beta\left(\xi_{\delta},\xi_{\rho}\right)=\eta_\beta\left(\left[\xi_{\delta},\xi_{\rho}\right]\right)=0$.
Finally, if $X=\xi_\delta$, $Y\in\Gamma\left(\cal H\right)$, then $d\eta_\beta\left(\xi_{\delta},Y\right)=0$
by Corollary \ref{facile2}. We conclude that $d\eta_\beta=0$ and  $M$ is  a $3$-cosymplectic manifold.
\end{proof}

\section{Corrected energy of $3$-quasi Sasakian manifolds}

As a single application, we show that the canonical foliation
of $3$-quasi Sasakian manifolds is a minimum of the \emph{corrected energy}.
Recall that the corrected energy ${\cal D}(\cal V)$ of a $p$-dimensional distribution $\cal V$
on a compact oriented Riemannian manifold $(M^m,g)$ is defined as (cf. \cite{chacon})
\begin{equation}\label{corrected}
{\cal D}({\cal
V})=\int_{M}\left(\sum_{a=1}^{m}\|\nabla_{e_a}\xi\|^2+q(q-2)\|\overrightarrow{H}_{\cal
H}\|^2+p^2\|\overrightarrow{H}_{\cal V}\|^2\right)d\textrm{vol},
\end{equation}
where $\left\{e_1,\ldots,e_m\right\}$ is a local adapted frame with
$e_1,\ldots,e_p\in{\cal V}_x$ and $e_{p+1},\ldots,$
$e_{m=p+q}\in{\cal H}_x={\cal V}_x^\perp$, and
$\xi=e_1\wedge\cdots\wedge e_p$ is a $p$-vector which determines the
distribution $\cal V$ regarded as a section of the Grassmann bundle
$G\left(p,M^m\right)$ of oriented $p$-planes in the tangent spaces
of $M^m$. Finally $\overrightarrow{H}_{\cal H}$ and
$\overrightarrow{H}_{\cal V}$ are the mean curvatures of the
distributions $\cal H$ and $\cal V$ given by
\begin{equation}
    \overrightarrow{H}_{\cal
H}=\sum_{\alpha=1}^{p}\left(\frac{1}{q}\sum_{i=p+1}^{n}h_{ii}^{\alpha}\right)e_\alpha,
\quad \quad \overrightarrow{H}_{\cal
V}=\sum_{i=p+1}^{n}\left(\frac{1}{p}\sum_{\alpha=1}^{p}h_{\alpha\alpha}^{i}\right)e_i,
    \label{HH}
\end{equation}
where the functions
$h_{ij}^{\alpha}=-g\left(\nabla_{e_i}e_\alpha,e_j\right)$ define the
second fundamental form of the horizontal distribution in the
direction $e_\alpha$, and
$h_{\alpha\beta}^{i}=-g\left(\nabla_{e_\alpha}e_i,e_\beta\right)$
define the second fundamental form of the vertical distribution in
the direction $e_i$. The norm of $\nabla \xi$ is given by
\begin{equation}\label{nablaxi}
\sum_a \|\nabla_{e_a}\xi\|^2=\sum_{i,j,\alpha}(h_{ij}^\alpha)^2
+\sum_{i,\alpha,\beta}(h_{\alpha\beta}^i)^2.
\end{equation}
It is proved in \cite{chacon} that if $\cal V$ is integrable then
\begin{equation}\label{inequality}
{\cal D(\cal V)}\geq\int_{M}\sum_{i,\alpha}c_{i\alpha}d\textrm{vol},
\end{equation}
where $c_{i\alpha}$ is the sectional curvature of the plane spanned
by $e_i\in\cal H$ and $e_\alpha\in\cal V$.

\begin{theorem}
The canonical foliation $\cal V$ defined by the Reeb vector fields
$\xi_1,\xi_2,\xi_3$ of a compact oriented $3$-quasi-Sasakian manifold
$(M,\phi_\alpha,\xi_\alpha,\eta_\alpha,g)$ represents a
minimum of the corrected energy ${\cal D}({\cal V})$.
\end{theorem}
\begin{proof}
We show that (\ref{inequality}) is in fact an equality. Fix an
adapted basis of type $\{e_1,\ldots,e_{4n},\xi_1,\xi_2,\xi_3\}$,
with $e_i\in\cal H$. Using \eqref{conseguenza1kanemaki} we have
$h_{ii}^\alpha= -g(\nabla_{e_i}\xi_\alpha, e_i)=-g(\phi_\alpha
A_\alpha e_i,e_i)$ which vanishes since $A_\alpha$ is symmetric and
commutes with $\phi_\alpha$. Moreover,
$h_{\alpha\beta}^i=-g(\nabla_{\xi_\alpha}e_i,\xi_\beta)
=g(\nabla_{\xi_\alpha}\xi_\beta,e_i)=0$ since the foliation $\cal V$
is totally geodesic. It follows that the mean curvatures
$\overrightarrow{H}_{\cal H}$ and $\overrightarrow{H}_{\cal V}$
vanish. Now we will express the norm of $\nabla \xi$ in terms of the
norms of the Reeb vector fields. We know that in a
$3$-quasi-Sasakian manifold $[\xi_\alpha,\xi_\beta]=c\xi_\gamma$,
where $c$ is a real number, possibly zero, and
$(\alpha,\beta,\gamma)$  is an even permutation  of $\{1,2,3\}$.
Applying the Koszul formula for the Levi-Civita connection, one
easily gets $\nabla_{\xi_\alpha}\xi_\beta=\frac{c}{2}\xi_\gamma$.
Then, we have
\begin{align*}
\|\nabla\xi_\alpha\|^2 &= \sum_i
g(\nabla_{e_i}\xi_\alpha,\nabla_{e_i}\xi_\alpha)
+g(\nabla_{\xi_\beta}\xi_\alpha,\nabla_{\xi_\beta}\xi_\alpha)
+g(\nabla_{\xi_\gamma}\xi_\alpha,\nabla_{\xi_\gamma}\xi_\alpha)\\
&=\sum_i g(\phi_\alpha A_\alpha e_i,\phi_\alpha A_\alpha
e_i)+\frac{c^2}{2}.
\end{align*}
On the other hand, applying (\ref{nablaxi}) and
$h_{\alpha\beta}^i=0$, we get
$$\sum_a
\|\nabla_{e_a}\xi\|^2=\sum_{i,j,\alpha}g(\nabla_{e_i}\xi_\alpha,e_j)^2=
\sum_{i,j,\alpha}g(\phi_\alpha A_\alpha
e_i,e_j)^2=\sum_{i,\alpha}g(\phi_\alpha A_\alpha e_i,\phi_\alpha
A_\alpha e_i),$$ since $A_\alpha {\cal H}\subset {\cal H}$ by Corollary \ref{preserva}. It follows that
$$\sum_a
\|\nabla_{e_a}\xi\|^2=\sum_{\alpha=1}^3
\|\nabla\xi_\alpha\|^2-\frac{3}{2}c^2.$$ Therefore the expression of
the corrected energy of $\cal V$, given by \eqref{corrected},
reduces to
\begin{equation*}
{\cal D}({\cal
V})=\int_{M}\left(\sum_{\alpha=1}^{3}\|\nabla\xi_\alpha\|^2-\frac{3}{2}c^2\right)
d\textrm{vol}.
\end{equation*}
Now, a direct computation shows that the sectional curvature of the
plane spanned by $\xi_\alpha$ and $\xi_\beta$ is
$K(\xi_\alpha,\xi_\beta)=\frac{c^2}{4}$. Hence, by
\eqref{olszakricci} we get
\begin{align*}
\sum_{i=1}^{4n}\sum_{\alpha=1}^{3}c_{i\alpha}&=
\sum_{\alpha=1}^{3}\sum_{i=1}^{4n}K\left(e_i,\xi_{\alpha}\right)\\
&=\sum_{\alpha=1}^{3}\left(\textrm{Ric}\left(\xi_\alpha\right)
-K(\xi_\alpha,\xi_\beta)-K(\xi_\alpha,\xi_\gamma)\right)\\
&=\sum_{\alpha=1}^{3}\|\nabla\xi_\alpha\|^2-\frac{3}{2}c^2,
\end{align*}
from which the assertion follows.
\end{proof}

\bigskip
Department of Mathematics,  University of Bari \\
  Via Edoardo Orabona, 4\\
  I-70125 Bari (Italy)\\
\\
\texttt{cappelletti@dm.uniba.it}\\
\texttt{antondenicola@gmail.com}\\
\texttt{dileo@dm.uniba.it}
\end{document}